\documentclass{gtart_h}

\def\ifplaintex{\expandafter\ifx\csname documentclass\endcsname\relax}

\def\gtp{{\mathsurround=0pt\it $\cal G\mskip-2mu$eometry \&\ 
$\cal T\!\!$opology $\cal P\!$ublications}}  

\def\recd{{\small Received:\qua\receiveddate\ifx\reviseddate\relax
\else\qquad Revised:\qua\reviseddate\fi\par}} 

\def\lognumber#1{\def\thelognumber{#1}}
\def\volumenumber#1{\def\thevolumenumber{#1}}
\def\volumeyear#1{\def\thevolumeyear{#1}}
\def\papernumber#1{\def\thepapernumber{#1}}
\def\pagenumbers#1#2{\def\startpage{#1}\def\finishpage{#2}}
\def\published#1{\def\publishdate{#1}}

\def\received#1{\def\receiveddate{#1}}
\def\revised#1{\def\reviseddate{#1}}
\def\accepted#1{\def\accepteddate{#1}}

\def\asciiurl#1{\def\theasciiurl{#1}}

\def\asciikeywords#1{\def\theasciikeywords{#1}}

\let\\\par\let\thelognumber\relax\let\thevolumenumber\relax
\let\thepapernumber\relax\let\thevolumeyear\relax\let\startpage\relax
\let\finishpage\relax\let\publishdate\relax\let\receiveddate\relax
\let\reviseddate\relax\let\accepteddate\relax\let\theasciititle\relax
\let\theasciiauthors\relax
\let\theasciiabstract\relax\let\theasciikeywords\relax

\let\theasciiemail\relax
\let\theasciiurl\relax

\ifplaintex
\font\logobig=cmssbx10 scaled 3836
\font\logomed=cmssbx10 scaled 2557
\else
\font\logobig=cmssbx10 scaled 4200
\font\logomed=cmssbx10 scaled 2800
\fi

\long\def\makeagttitle{   
\count0=\startpage
\agt\hfill      
\hbox to 45truept{\vbox to 0pt{\vglue -13truept{\logomed A\kern -.37em{\logobig 
T}\kern -.38em G}\vss}\hss}
\break
{\small Volume \thevolumenumber\ (\thevolumeyear)
\startpage--\finishpage\nl
Published: \publishdate}

\vglue .25truein

{\parskip=0pt\leftskip 0pt plus
1fil\def\\{\par\smallskip}{\Large\bf\thetitle}\par\medskip} \vglue
0.05truein

{\parskip=0pt\leftskip 0pt plus 1fil\def\\{\par}{\sc\theauthors}
\par\medskip}%
 
\vglue 0.03truein 

{\small\leftskip 25truept\rightskip 25truept{\bf Abstract}\stdspace\theabstract

{\bf AMS Classification}\stdspace\theprimaryclass
\ifx\thesecondaryclass\relax\else; \thesecondaryclass\fi\par
{\bf Keywords}\stdspace \thekeywords\par}\vglue 7truept

}   

\ifplaintex
\font\phead=cmsl9 scaled 950
\font\pnum=cmbx10 scaled 913
\font\pfoot=cmsl9 scaled 950
\headline{\vbox to 0pt{\vskip -4.5mm\line{\small\phead\ifnum
\count0=\startpage ISSN 1472-2739 (on-line) 1472-2747 (printed)
\hfill {\pnum\folio}\else\ifodd\count0\def\\{ }%
\ifx\theshorttitle\relax\thetitle\else\theshorttitle\fi\hfill{\pnum\folio}
\else\def\\{ and }{\pnum\folio}\hfill\ifx\theshortauthors\relax\theauthors
\else\theshortauthors\fi\fi\fi}\vss}}
\footline{\vbox to 0pt{\vglue 0mm\line{\small\pfoot\ifnum\count0=\startpage
\copyright\ \gtp\hfill\else
\agt, Volume \thevolumenumber\ (\thevolumeyear)\hfill\fi}\vss}}
\else
\font=cmsl9 scaled 1050
\font\lnum=cmbx10 
\font=cmsl9 scaled 1050
\makeatletter
\def\@oddhead{{\small\lhead\ifnum\count0=\startpage ISSN 1472-2739 
(on-line) 1472-2747 (printed)\hfill {\lnum\number\count0}\else\ifodd\count0
\def\\{ }\ifx\theshorttitle\relax \thetitle \else\theshorttitle\fi\hfill
{\lnum\number\count0}\else\def\\{ and }{\lnum\number\count0}
\hfill\ifx\theshortauthors\relax 
\theauthors\else\theshortauthors\fi\fi\fi}}\def\@evenhead{\@oddhead}
\def\@oddfoot{\small\lfoot\ifnum\count0=\startpage\copyright\ \gtp\hfill\else
\agt, Volume \thevolumenumber\ (\thevolumeyear)\hfill\fi}
\def\@evenfoot{\@oddfoot}
\makeatother
\fi
\let\maketitlepage\makeagttitle
\let\makeshorttitle\maketitlepage
\let\maketitle\maketitlepage

\newwrite\gtoutfile
\long\gdef\makeheadfile{  
{\def\\{, }\def\s{ }
\immediate\openout\gtoutfile head.xxx
\immediate\write\gtoutfile{Proxy-for: \ifx\theasciiauthors\relax
\theauthors\else\theasciiauthors\fi\s<\ifx\theasciiemail\relax\theemail\else\theasciiemail\fi>}
\immediate\write\gtoutfile{\noexpand\\}
\immediate\write\gtoutfile{Authors: \ifx\theasciiauthors\relax
\theauthors\else\theasciiauthors\fi}
{\def\\{ }\immediate\write\gtoutfile{Title: \ifx\theasciititle\relax
\thetitle\else\theasciititle\fi}}
\immediate\write\gtoutfile{Subj-class: GT or SG, GR etc}
\immediate\write\gtoutfile{MSC-class: \theprimaryclass\ifx\thesecondaryclass\relax\else, \thesecondaryclass\fi}
\immediate\write\gtoutfile{Journal-ref: Algebr. Geom. Topol. \thevolumenumber\s
(\thevolumeyear) \startpage-\finishpage}
\immediate\write\gtoutfile{Comments: Published by Algebraic and
Geometric Topology at}
\immediate\write\gtoutfile{\s\s\s  http://www.maths.warwick.ac.uk/agt/AGTVol\thevolumenumber/agt-\thevolumenumber-\thepapernumber.abs.html}
\immediate\write\gtoutfile{\noexpand\\}
\immediate\write\gtoutfile{}
\ifx\theasciiabstract\relax
\immediate\write\gtoutfile{\theabstract}\else
\immediate\write\gtoutfile{\theasciiabstract}\fi
\immediate\write\gtoutfile{}
\immediate\write\gtoutfile{\noexpand\\}
\immediate\write\gtoutfile{}
\immediate\closeout\gtoutfile}}  

\def\maketitlepage{\makeagttitle\makeheadfile}
\let\makeshorttitle\maketitlepage
\let\maketitle\maketitlepage

\lognumber{40}
\volumenumber{4}
\volumeyear{2004}
\papernumber{40}
\pagenumbers{935}{942}
\received{27 March 2001} 
\revised{27 September 2004}
\accepted{28 September 2004}
\published{13 October 2004}

\usepackage{epsfig,amsmath,amssymb}

\theoremstyle{plain}

\newtheorem{theorem}{Theorem}
\newtheorem{proposition}{Proposition}[section]

\newtheorem{corollary}[proposition]{Corollary}

\theoremstyle{definition}

\theoremstyle{remark}

\newcommand{\psdraw}[2]
         {\begin{array}{c} \hspace{-1.3mm}
        \raisebox{-4pt}{\epsfig{figure=draws/#1.eps,width=#2}}
        \hspace{-1.9mm}\end{array}}

\newlength{\standardunitlength}
\def\lbl#1{\label{#1}}

\def\BZ{\mathbb Z}
\def\BQ{\mathbb Q}

\def\A{\mathcal A}

\def\D{\Delta}
\def\K{\mathcal K}

\def\La{\Lambda}
\def\l{\lambda}

\def\ga{\gamma}

\def\fti{finite type invariant}

\def\o1o{\underset{1}\ast}
\def\x1y{\underset{x1y}\ast}
\def\y1x{\underset{y1x}\ast}

\def\ti{\widetilde}

\def\e{\varepsilon}

\def\AS{\mathrm{AS}}
\def\IHX{\mathrm{IHX}}
\def\clover{clover}  

\def\ylink{$\mathrm{Y}$-link}

\def\Wh{\mathrm{Wh}^{\epsilon}}

\def\lgt{\mathrm{lk}^{\ga}_{\ti M}}
\def\tKb{\tau_{\mathfrak{p}}}
\def\Zrat{Z^{\mathrm{rat}}}
\def\sub{\subset}
\def\eqI{\equiv^I}
\def\eqY{\equiv^Y}
\def\eql{\equiv^l}
\def\Th{\Theta}

\begin{document}


\title[Whitehead doubling persists]{Whitehead doubling persists}

\author{Stavros Garoufalidis}
\address{School of Mathemtaics, Georgia Institute of 
Technology\\Atlanta, GA 30332-0160, USA. }
\email{stavros@maths.gatech.edu}
\urladdr{http://www.math.gatech.edu/~stavros } 
\asciiurl{http://www.math.gatech.edu/ stavros }

\asciikeywords{Whitehead double, loop filtration,
Goussarov-Habiro, clovers, claspers, Kontsevich integral}

\begin{abstract}
The operation of (untwisted) Whitehead doubling trivializes the Alexander
module of a knot (and consequently, all known
abelian invariants), and converts knots to topologically slice ones. In this
note we show that Whitehead doubling does not trivialize the rational function
that equals to the 2-loop part of the Kontsevich integral.
\end{abstract}

\maketitle


\section{Introduction}
\lbl{sec.intro}

\subsection{History}
\lbl{sub.history}

The {\em colored Jones function} $J$ of a knot $K$ is a 2-parameter formal 
power series
$$
J(K)(h,\l)=\sum_{n,m=0}^{\infty} a_{n,m}(K) h^n \l^m
$$
with remarkable periodicity properties. By its definition, if $\l=d$ is a 
natural number, $J(h,d)$ coincides with the colored Jones polynomial of 
the knot (using the $(d+1)$-dimensional irreducible representation of 
$\mathfrak{sl}_2$) and thus 
$$
J(h,d) \in \BZ[e^{\pm h}].
$$
We can think of this as a {\em periodicity property} (i.e., a set of 
recursion relations) for the coefficients $a_{n,m}$ of $J$. This is
an obvious periodicity property.

We now come to describe some {\em hidden periodicity} of $J$. Each 
coefficient $a_{n,m}$ is a Vassiliev invariant of degree $n$ and vanishes 
if $m > n$. Thus, we can rearrange $J$ as a sum of subdiagonal terms
$$
J(h,\l)=\sum_{k=0}^\infty h^k Q_{J,k}(h \l)
$$
where
$$
Q_{J,k}(s)=\sum_{m=0}^\infty a_{k+m,m} s^m.
$$ 
The MMR conjecture, shown in \cite{BG}, states that $Q_{J,0}$ is a 
reparametrization of the {\em Alexander polynomial} of the knot. This
translates to a hidden periodicity of $J$.

Rozansky, in his study of the colored Jones function conjectured that
for every $n \geq 1$, $Q_n$ is a reparametrization of a rational function,
whose denominator is a power of the Alexander polynomial. Using a variety
of quantum field theory techniques and an appropriate expansion of the
$R$-martix, in a difficult paper Rozansky proved the above mentioned
conjecture, \cite{Ro1}.

It is well-known that the colored Jones function is an image of the
graph-valued universal Vassiliev invariant, the so-called {\em Kontsevich
integral}. The Kontsevich integral has a subdiagonal expansion, and 
Rozansky further conjectured that each term of the expansion of the Kontsevich
integral should be given by a series of trivalent graphs with rational 
functions attached to their edges. This is often called the {\em Rationality
Conjecture}.

A weak form of the Rationality Conjecture was quickly proven by Kricker 
\cite{Kr}. A stong form followed by joint work with Kricker \cite{GK}, where 
a rational noncommutative invariant $\Zrat$ of knots was constructed.

Although the construction of the $\Zrat$ invariant is rather involved, there
are several axiomatic properties which make it easier to understand parts
of the $\Zrat$ invariant in terms of geometric invariants of knots.

For example, consider the move that replaces a knot $K$ in an integer
homology sphere by a knot $K'$ obtained by surgery on a hullhomologous
clasper in the complement of $K$. This {\em null-move} was introduced in 
\cite{GR}. For a reference on claspers, see \cite{Gu1, Gu2, Ha} and also
\cite{GGP}.
As was shown in \cite{GK}, each term of the $\Zrat$ is a \fti 
with respect to the null move. Moreover, $0$-equivalence (under the null
move) coincides with $S$-equivalence. Below, we will evaluate $\Zrat$ on a set
of $S$-equivalent knots. In this case, the matrix part of $\Zrat$ is fixed
and its graph part takes values in a graded vector space, see 
\cite{GK}. Let $Q_n=\Zrat_n$ denote the degree $n$ part of 
$\Zrat$. 

$Q_n$ takes values in a vector space over $\BQ$ generated
by trivalent graphs with $2n$ trivalent vertices. The graphs have rational
functions in $t$ assigned to their edges, and satisfy certain linear relations
explained in \cite{GR} and \cite{GK}.
For the statement and the proof of Theorem \ref{thm.1} below, we only need 
to know that $Q_n$ takes value in an {\em abelian} group; and any vector
space over $\BQ$ is an abelian group.

\subsection{Statement of the results}
\lbl{sub.results}

Our first result concerns the change of $Q_n$ under a modification of a knot.
All knots will be oriented and, unless otherwise mentioned, $0$-framed.
Consider a knot $K_0$ which intersects a ball $B \sub S^3$ in two unknotted
arcs with opposite orientation. This {\em pattern} $\mathfrak{p}=(K_0,B)$ 
gives rise to a map:
$$
\tKb: \text{Knots} \to \text{Knots}
$$
which sends $K$ to the result of replacing $B \cap K_0$ with a 2-parallel
of a $1$-tangle version of $K$, with $0$-framing.

This move on knots can be described in terms of surgery in the ambient space 
and has a long history. It was already used in the sixties (under the term,
surgical modification, see for example Levine \cite{L1}) 
to prove realization theorems for algebraic obstructions. 
Modification of knots was also used in the seventies by Casson-Gordon in 
their secondary obstruction invariants. Later on, the theme was taken
by many authors including Gilmer, Livingston, and Cochran-Orr-Teichner (the
latter, use the biological term: infection). In the world
of finite type invariants, a systematic study of geometric surgery on 
clovers or claspers was initiated by Goussarov and Habiro \cite{Gu1,Ha}.

\begin{theorem}
\lbl{thm.1}
{\rm(i)}\qua For fixed $K_0$ and $n \geq 1$, the map 
$$
\phi_n: K \to Q_n(\tKb(K))
$$
is a \fti \ of 0-framed knots of type $2n$, whose 
degree $2n$ part lies in the algebra of Alexander-Conway coefficients.
\newline
{\rm(ii)}\qua  For $n=1$, we have that
$$Q_1(\tKb(K))=Q_1(\tKb(\mathrm{unknot}))+ c_{\mathfrak{p}} \cdot a(K)
$$
where $a(K)$ is a nontrivial Vassiliev invariant of degree $2$
(such as the second derivative of the Alexander polynomial) and 
$c_{\mathfrak{p}}$ is a constant that depends on the pattern.
\end{theorem}

As an application, consider the following pattern:
$$
\psdraw{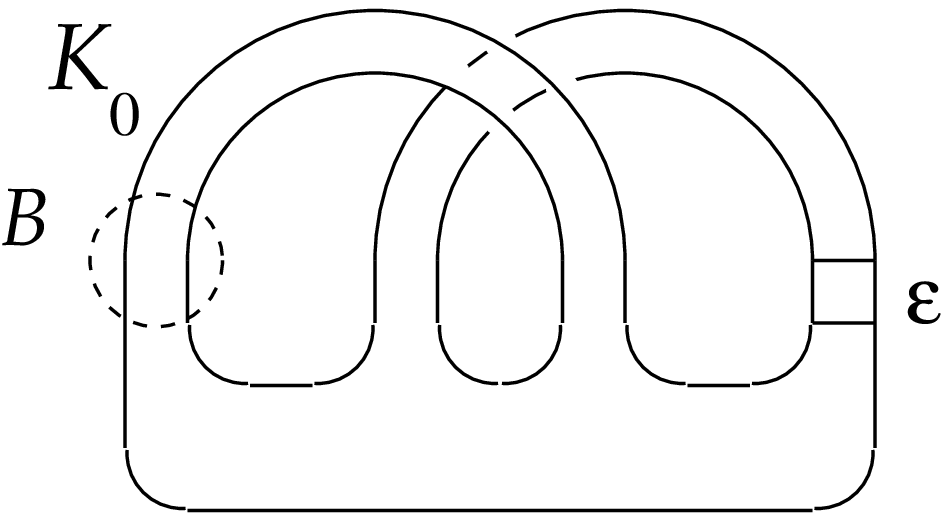}{1.5in}
$$
(where the $\e=\pm 1$ indicates a full twist, depending on the sign
of the clasp).
Then, 
$$
\tKb(K)=\Wh(K)
$$
is the {\em untwisted Whitehead double} of $K$ with either 
clasp, i.e., the satellite of $K$ with respect to the pattern:
$$\psdraw{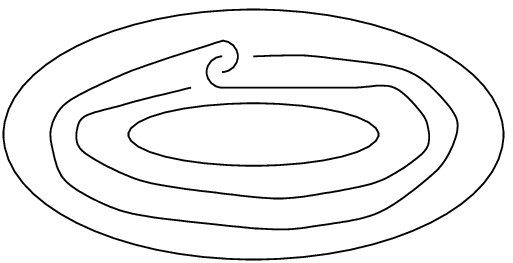}{1.5in}.$$ 

Here is a sample corollary. Since it involves explicit computations
with the $Q_1$ invariant evaluated on knots with trivial Alexander 
polynomial, we need to recall that
$$
Q_1: \text{Knots with trivial Alexander polynomial} \to
\A_2(\La)
$$
where
$$
\A_2(\La)=\otimes^3 \BQ[t^{\pm 1}]/((f,g,h)=(tf,tg,th), 
\text{Aut}(\Th))
$$ 
is the quotient of the abelian group 
$\otimes^3 \BQ[t^{\pm 1}]$ modulo the subgroup generated by 
$a \otimes b \otimes c=ta \otimes tb \otimes tc$ 
and the subgroup generated by $a \otimes b \otimes c= ga \otimes gb 
\otimes gc$ for $g \in \text{Aut}(\Th)$.
Here, $\text{Aut}(\Th)=\text{Sym}_3 \times \text{Sym}_2$ which
acts on $\otimes^3 \BQ[t^{\pm 1}]$ by permutation on the 
3 factors of the tensor product and by simultaneous replacement of $t$
(in all factors of the tensor product) by $t^{-1}$.

In order to make contact with \cite{GR} and \cite[Sec.5.3]{GK}, and to explain 
the origin of the $\text{Aut}(\Th)$ symmetry, we point out a graphical 
interpretation  of elements of $\A_2$ by trivalent graphs with oriented edges 
and beads (that is, elements of $\BQ[t^{\pm 1}]$) on their edges:
$$
a \otimes b \otimes c \leftrightarrow 
\psdraw{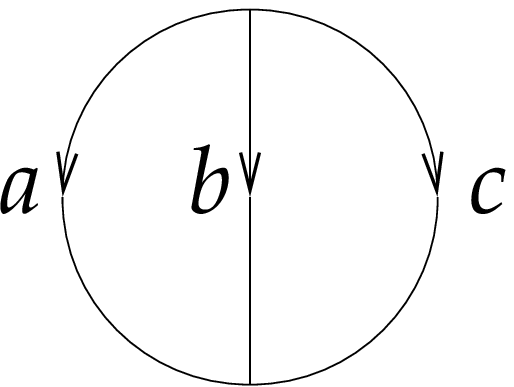}{1in} .
$$

\begin{corollary}
\lbl{cor.1}
We have:
$$
Q_1(\Wh(K))= \e a(K) \cdot  1 \otimes 1 \otimes (t+t^{-1}-2)
$$
where $a(K)=1/4\frac{d^2}{dh^2}|_{h=0}\D(K)(e^h) \in \BZ$, and $\e=\pm 1$ is
the sign of the clasp.
\end{corollary}

Whitehead doubling trivializes the Alexander module
as well as the more sophisticated topological slicing obstructions
of Casson-Gordon and Cochran-Orr-Teichner.
On the other hand, $Q_1$ remembers better the knot that is about to be
Whitehead doubled.

\section{Proofs}

\begin{proof}[Proof of Theorem \ref{thm.1}]$\phantom{9}$

Let us recall three different moves on the set $\K$ of 0-framed knots in 
$S^3$: 
\begin{itemize}
\item
Changing a crossing, i.e., doing an $I$-modification in the language of
\cite{Gu2} and \cite{Ha}.
\item
Doing a $\Delta$-move, i.e., doing a $Y$-modification in the language of
\cite{MN}.
\item
Doing a null-move, in the language of \cite{GR}.
\end{itemize}
These three moves lead in the usual way to three notions of \fti s and 
corresponding notions of 
$n$-{\em equivalence}, denoted by $\eqI_n, \eqY_n$ and  $\eql_n$ 
respectively. Note that $K \equiv_n K'$ implies that $f(K)=f(K')$
for all invariants $f$ of type $n$. It is a folk result (easily proven
by the results of Goussarov and Habiro) that $K \eqI_{n+1} K'$ iff
$K \eqY_n K'$. Furthermore, it is easy to see that if $K \eqY_n K'$,
then $\tKb(K) \eql_n \tKb(K')$. Further, in \cite{GR} and \cite{GK},
itt was shown that $Q_n$ is a  \fti \ of type
$2n$ with respect to the null move, where $Q_n$ takes values in an appropriate 
$\BQ$-vector space.

This discussion implies the following conclusion, for every fixed $n$. 
$$
K \eqI_{2n+1} K' \implies K \eqY_{2n} K' \implies \tKb(K) \eql_{2n} \tKb(K')
\implies  \phi_n(K)=\phi_n(K').
$$
Thus, $\phi_n: K \mapsto Q_n(\tKb(K))$ 
is an additive (under connected sum) function on $\K/\K^I_{2n+1}$ 
(here $\K^I_n$ denotes the
set of $n$-trivial knots with respect to the $I$-move). By a theorem of 
Goussarov and Habiro it follows that $\phi_n$ is a Vassiliev invariant of 
degree at most $2n+1$. We claim that $\phi_n$ is $\BQ$-valued of Vassiliev 
degree $2n$. Indeed, $\K^I_{2n}/\K^I_{2n+1} \otimes \BQ$ is a vector space 
spanned by uni-trivalent graphs $G$ with $2(2n+1)$ vertices such that every 
component of $G$ contains a trivalent vertex, modulo the $\AS$ and $\IHX$ 
relations. Using the $\AS$ relation, we can assume that distinct univalent 
vertices of $G$ are joined to distinct trivalent vertices. It follows that
the number of trivalent vertices of $G$ is at least equal to the number of 
univalent vertices, and thus that $G$ has at least $2n+1$ trivalent vertices, 
which gives rise to a null move of degree $2n+1$, on which $Q_n$ vanishes.

This implies that $\phi_n$ is of Vassiliev degree $2n$. Furthermore, if
$G$ is a unitrivalent graph of degree $2n$ with more than $2n$ trivalent 
vertices, then $\phi_n(G)=0$. The remaining graphs of degree $2n$ are
a disjoint union of {\em wheels} i.e., diagrams like 
$\psdraw{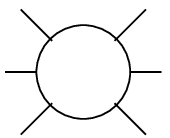}{0.2in}$. This, together with a result of 
Kricker-Spence-Aitchinson \cite{KSA} implies that
the degree $2n$ part of $Q_n$ lies in the algebra of Alexander-Conway
coefficients. 

This concludes the proof of the first part Theorem \ref{thm.1}.
The second part follows easily since, up to a multiple, there is a unique
$\BQ$-valued Vassiliev invariant of $I$-degree $2$.
\end{proof}

\begin{proof}[Proof of Corollary \ref{cor.1}]$\phantom{9}$

Notice that $\Wh(\text{unknot})=\text{unknot}$ and that $Q_1(\text{unknot})
=0$. Thus, $\phi_1(K)=c a(K)$, where $a(K)=\frac{d^2}{dt^2}|_{t=1}
\D(K)(t)$ 
and $c$ is a constant. In order to figure out $c$, we need a computation.
Let $G=(G_1,G_2)$ be a wheel with two legs attached to an unknot $K$.
Consider the \ylink \ of degree $2$ (also denoted by $G$) in the complement
of $\Wh(K)$. Let $L_{ij}$ for $i=1,2$ and $j=1,2,3$ denote the six leaves
of $G$ labeled as in Figure \ref{G}. 
 
\begin{figure}[htpb]
$$\psdraw{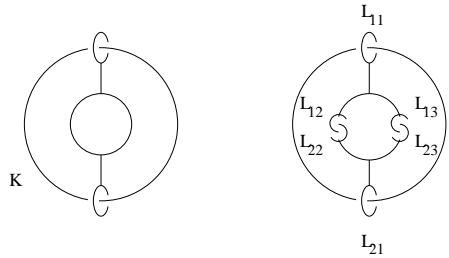}{2.3in}$$
\caption{On the left, a wheel with two legs. On the right, the associated
\clover \ of degree $2$}
\lbl{G}
\end{figure}

It follows from \cite{GR} (see also \cite{GK}) that
$
Q_1([(S^3,\Wh(K)),G]) 
$
is the result of complete contractions of all pairs of leaves of 
$G_1 \cup G_2$, where we label the edges by equivariant linking numbers
of the leaves. This gives rise to a linear combination of trivalent
graphs with two trivalent vertices and edges decorated by elements of 
$\BZ[t,t^{-1}]$. 

The equivariant linking function $\lgt$ (see \cite{GR}) of the leaves 
is given as follows:
$$
\lgt(L_{1i},L_{2j})=
\begin{cases}
\delta_{i,j} & \text{ for } i=2,3 \\
\e \delta_{1,j}(t+t^{-1}-2) & \text{ for } i=1,
\end{cases}
$$
where $\delta_{a,b}=1$ (resp. $0$) if $a=b$ (resp. $a \neq b$),
and the sign $\e=\pm 1$ depends on the sign of the clasp of the 
Whitehead double.

The computation of the equivariant linking function is best seen by drawing 
the universal abelian cover of $\Wh(K)$ and lifting $G$ to it.  

Alternatively, one may use the genus $1$ Seifert surface of $\Wh(K)$ drawn 
above and identify the equivariant linking numbers in question with the 
equivariant linking numbers of links formed by 
meridians dual to the bands. 

In \cite[Prop.14.3]{L2} Levine computes the
matrix $B=(\lgt(m_i,m_j))_{i,j}$ of equivariant linking numbers of meridians 
$m_i$ dual to the bands of a Seifert surface by:
$$
B=(t-1)(tA-A^T)^{-1}
$$
where $A$ is the Seifert matrix with respect to a basis consisting of bands,
and $A^T$ is the transpose of $A$. In our case, the Seifert matrix is
$$
A= \left( \begin{matrix}
0 & 1 \\ 0 & \e
\end{matrix} \right)
$$
and the corresponding matrix $B$ is
$$
A= \left( \begin{matrix}
\e (t+t^{-1}-2) & 1-t \\ 1-t^{-1} & 0
\end{matrix} \right) .
$$

On the other hand, we have that $a([K,G])=1$, since
$$\D(K)(e^h)=\exp(-2\sum_{n} a_{2n}(K) h^{2n}),$$ where $a_{2n}(K)$ is
the coefficient of the degree $2n$ wheel $w_{2n}$ in the logarithm
of the Kontsevich integral, \cite{KSA}.
The result follows.
\end{proof}

\rk{Acknowledgments}The author is partially supported by NSF grant
        DMS-98-00703 and by an Israel-US BSF grant.

\Addresses\recd

\end{document}